\documentclass{amsart}
\usepackage{calc}

\usepackage{amsmath, amsfonts, amsthm, amssymb, color,  graphicx, latexsym, mathrsfs, cite}
\usepackage{comment,stmaryrd}
\usepackage{ifpdf}
\ifpdf \usepackage[colorlinks=true, citecolor=blue, linkcolor=blue, urlcolor=blue]{hyperref} \fi

\newcommand{\cal}{\mathcal}

\newtheorem{formula}{}[section]
\newtheorem{definition}[formula]{Definition}
\newtheorem{corollary}[formula]{Corollary}
\newtheorem{remark}[formula]{Remark}
\newtheorem{lemma}[formula]{Lemma}
\newtheorem{theorem}[formula]{Theorem}

\def\thrm{\begin{theorem}}
\def\thrml#1{\begin{theorem}\label{#1}}
\def\ethrm{\end{theorem}}
\def\rmrk{\begin{remark}}
\def\rmrkl#1{\begin{remark}\label{#1}}
\def\ermrk{\end{remark}}
\def\dfntn{\begin{definition}}
\def\dfntnl#1{\begin{definition}\label{#1}}
\def\edfntn{\end{definition}}
\def\nmrt{\begin{enumerate}}
\def\enmrt{\end{enumerate}}
\def\tm#1{\item[{\rm (#1)}]}
\def\qtnl#1{\begin{equation}\label{#1}}
\def\eqtn{\end{equation}}
\def\lmm{\begin{lemma}}
\def\lmml#1{\begin{lemma}\label{#1}}
\def\elmm{\end{lemma}}
\def\crllr{\begin{corollary}}
\def\crllrl#1{\begin{corollary}\label{#1}}
\def\ecrllr{\end{corollary}}
\def\css{\begin{cases}}
\def\ecss{\end{cases}}

\def\proof{\noindent{\bf Proof}.\ }

\def\cX{{\cal X}}

\def\mN{{\mathbb N}}

\DeclareMathOperator{\aut}{Aut}

\DeclareMathOperator{\GCD}{GCD}

\DeclareMathOperator{\id}{id}
\DeclareMathOperator{\im}{im}

\DeclareMathOperator{\orb}{Orb}

\DeclareMathOperator{\syl}{Syl}
\DeclareMathOperator{\sym}{Sym}

\DeclareMathOperator{\zel}{Zel}

\def\eprf{\hfill$\square$}

\def\qaq{\quad\text{and}\quad}
\def\qoq{\quad\text{or}\quad}
\def\ov{\overline}
\def\phmb#1{{\phantom{x}\hspace{-2mm}^{#1}}}

\begin{document}
\title{On $2$-closed abelian permutation groups}
\author{Dmitry Churikov}
\author{Ilia Ponomarenko}
\thanks{The work is supported by Mathematical Center in Akademgorodok, the agreement with Ministry of Science and High Education of the Russian Federation number  075-15-2019-1613 and partially supported by the Russian Foundation for Basic Research (Grant 18-01-00752). }
\date{}

\begin{abstract}
A permutation group $G\le\sym(\Omega)$ is said to be $2$-closed if no  group $H$ such that $G<H\le\sym(\Omega)$ has the same orbits on $\Omega\times\Omega$ as~$G$. A simple and efficient inductive criterion for the $2$-closedness is established  for abelian permutation groups with cyclic transitive constituents.
\end{abstract}

\maketitle

\section{Introduction}

The concept of $m$-closed permutation groups, $m\in\mN$, was introduced by H.~Wie\-landt~\cite{WieInvRel} in the framework of the method of invariant relations, developed him to study group actions. He proved that the $m$-closed groups are exactly the automorphism groups of the relational structures formed by families of $m$-ary relations on the same  set. The $1$-closed groups are just the direct product of symmetric groups. In the present paper, we are interested in the $2$-closed groups.\medskip

Let $\Omega$ be a finite set and $G\le \sym(\Omega)$. Denote by $\orb_2(G)$ the set of all orbits of the induced action of~$G$ on~$\Omega\times\Omega$. The {\it $2$-closure} of the permutation group~$G$ is defined to be the largest subgroup $\ov G=G^{(2)}$  in $\sym(\Omega)$ that have the same $2$-orbits as $G$, i.e.,
\qtnl{260920b}
\orb_2(G)=\orb_2(\ov G).
\eqtn
The group $G$ is said to be {\it $2$-closed} if $G=\ov G$. The set $\orb_2(G)$ forms a relational structure on $\Omega$, which is called a coherent configuration~\cite{CCBook}. Thus the $2$-closed groups are the automorphism groups of coherent configurations. This fact makes the concept of $2$-closed groups especially important for algebraic combinatorics~\cite{EvdokimovP2009}.\medskip

Based on the definition alone, it is difficult to determine whether a given permutation group $G$ is $2$-closed. In some special cases (mostly for transitive groups the degree of which is the product of small number of primes), good criteria for the $2$-closedness are obtained via the full classification of  permutation groups of given degree, see, e.g.,~\cite{DK2009}. To the best of our knowledge, efficient algorithms  recognizing $2$-closed groups are known only for a few cases, e.g., for odd order groups~\cite{EvdokimovP2001} and for supersolvable groups~\cite{VasP2020}.\medskip

The present paper is motivated by the lack of good criteria for the $2$-closedness even for abelian groups. A criterion found in \cite[Theorem~6.1]{GK2019a} does not seem quite satisfactory, because requires for a group~$G$ to inspect all permutations in the direct product of the constituents $G^\Delta$, $\Delta\in\orb(G)$. Our first result reduces the $2$-closedness question to the case of~$p$-groups, and generalize a Wielandt observation that the classes of $p$-groups and abelian groups are invariant with respect to taking the $2$-closure.

\thrml{090820a}
Let $G$ be a nilpotent permutation group. Then $\ov G$ is nilpotent. Moreover,
$$
\ov G=\prod_{P\in\syl(G)}\ov P.
$$
\ethrm

\crllrl{020619a}
A nilpotent permutation group $G$ is $2$-closed if and only if every Sylow subgroup of $G$ is $2$-closed.
\ecrllr

Probably, the first attempt to classify the $2$-closed abelian groups was made by A.~Zelikovskii in~\cite[Corollary~5]{Z89}.\footnote{The cited statement was formulated for arbitrary intransitive groups with $2$-closed transitive constituents.} However, the characterization appears to be wrong and infinitely many counterexamples were found in~\cite{GK2018}.  To explain the gap in ~Zelikovskii's argument, let $G\le\sym(\Omega)$. One can associate with $G$  a permutation group on $\Omega$, defined as follows:
\qtnl{270820a}
\zel(G)=\prod_{\Delta} \bigcap_{\Delta'\ne\Delta} G_{\scriptscriptstyle\Delta'}^{\scriptscriptstyle\Delta},
\eqtn
where $\Delta$ and $\Delta'$ run over the orbits of $G$, and $G_{\scriptscriptstyle\Delta'}^{\scriptscriptstyle\Delta}=(G_{\scriptscriptstyle\Delta'})^{\scriptscriptstyle\Delta}$ is the restriction to~$\Delta$ of the pointwise stabilizer of~$\Delta'$ in~$G$. In this notation, the necessary and sufficient condition claimed by Zelikovskii for an abelian group $G$ to be $2$-closed is that
$$
\zel(G)^\Delta=(G_{\Omega\setminus\Delta})^\Delta\quad\ \,\text{for all}\ \,\Delta\in\orb(G).
$$
It is not hard to see  that this condition is equivalent to the inclusion $\zel(G)\le G$. Essentially,  the only gap in Zelikovskii's argument consists in the wrong statement that if $G$ is $2$-closed, then so is the permutation group induced by the action of~$G$ on the orbits of~$\zel(G)$. In this way we arrive at the second result of the present paper.

\thrml{200420b}
Let $G$ be an abelian permutation group and $Z=\zel(G)$. Then $G$ is $2$-closed if and only if $Z\leq G$ and $G^{\orb(Z)}$ is $2$-closed.
\ethrm

The ``only if'' part of  Theorem~\ref{200420b} can be illustrated by the two following examples in which we construct non-$2$-closed groups~$G$ and~$H$, respectively. In the first example, $\zel(G)\not\leq G$, whereas in the second one, $\zel(H)\leq H$ but  $H^{\orb(\zel(H))}$ is not $2$-closed.\medskip

{\bf Example 1.} Let $p$ be a prime, and let $G\le\sym(3p)$ be an elementary abelian group of order~$p^2$. The action of~$G$ is chosen so that (a) there are exactly three $G$-orbits, each of size~$p$, and (b)  for any two points $\alpha$ and $\beta$ belonging to different $G$-orbits, the stabilizers $G_\alpha$ and $G_\beta$\ are different subgroups of~$G$ of order~$p$  (if $\alpha$ and~$\beta$ belong to the same orbit, then, of course, $G_\alpha=G_\beta$). From Wielandt's dissection theorem~\cite[Theorem~6.5]{WieInvRel}, it follows that the group~$\ov G$ equals the direct product of its transitive constituents; since $\ov G$ is abelian and $G\le\ov G$, we conclude that $|\ov G|=p^3$. Thus, $G$ is not $2$-closed and $\ov G=\zel(G)$.\medskip

{\bf Example 2.} Let $\Omega_1$ and $\Omega_2$ be disjoint sets, and let $G_1\le\sym(\Omega_1)$ and $G_2\le\sym(\Omega_2)$ be two copies of the permutation group $G$ from Example~1. For any two  permutations $g_1\in G_1$ and $g_2\in G_2$, corresponding to the same element of $G$, we define a permutation $g\in\sym(\Omega_1\cup \Omega_2)$ such that $g^{\Omega_i}=g_i$, $i=1,2$. Let $H$ be the group of all these permutations~$g$; this group is isomorphic to $G$ as abstract group and have the same point stabilizers as~$G$ (again as abstract groups). Then $H^\Delta_{\Delta'}=1$ for any two orbits $\Delta$ and $\Delta'$ corresponding to the same orbit of~$G$. It follows that $H$ is not $2$-closed (see Lemma~\ref{220920a}) and $\zel(H)=1$.\medskip

Theorem~\ref{200420b} reduces the question on the $2$-closedness of abelian group~$G$ to the case when the group $\zel(G)$ is trivial. This case can really happen even if $G$ is nontrivial, see Example~$2$. In a large class of abelian permutation groups, one can continue the reduction by ``removing'' {\it unessential orbits} defined as follows. An orbit $\Delta$ of a group $G\le\sym(\Omega)$ is said to be {\it unessential} if $G$ is $2$-closed if and only if $G^{\Omega\setminus\Delta}$ is $2$-closed.

\thrml{310820a}
Let $G$ be an intransitive group and $p$ a prime. Assume that every transitive  constituent of $G$ is a cyclic $p$-group. Then $\zel(G)$ is trivial only if every orbit of~$G$ is unessential.
\ethrm

Combining the obtained results we arrive at the following simple and efficient inductive criterion of the $2$-closedness for an (abelian) permutation group $G$ with cyclic transitive constituents. First, there is nothing to do if $G$ is transitive, because in this case $G$ is $2$-closed by~\cite[Example~5.13.]{WieInvRel}. Second, Theorem~\ref{090820a} reduces the problem to the case when~$G$ is a $p$-group. At this point a further reduction is needed. Namely, depending on whether the group $Z=\zel(G)$ is trivial, we continue with
$$
G^{\orb(Z)}\qoq G^{\Omega\setminus\Delta},
$$
where $\Delta$ is an arbitrarily chosen orbit of~$G$. The correctness is provided by  Theorems~\ref{200420b} and~\ref{310820a}, respectively. Each of these two reductions decreases the degree of a group and hence the test is finished after at most~$|\Omega|$ reductions.\medskip

The proofs of Theorems ~\ref{090820a}, \ref{200420b}, and~\ref{310820a} are presented in Sections~\ref{sec:Section3}, \ref{sec:Section4}, and~\ref{260920a}, respectively. The relevant notation and definitions are collected in Section~\ref{sec:Section2}.\medskip

The authors thank S.~Skresanov and A.~Vasil'ev for their suggestions for improving the text.

\section{Permutation groups}\label{sec:Section2}

Throughout the paper, $\Omega$ is a finite set. For a permutation group $G\le\sym(\Omega)$, we use the notation
$$
\orb(G)=\orb(G,\Omega)=\{\alpha^G:\ \alpha\in\Omega\}
$$
for the set of $G$-orbits $\alpha^G=\{\alpha^g:\ g\in G\}$. A set $\Delta\subseteq\Omega$ is said to be {\it $G$-invariant} if $\Delta^g=\Delta$ for all $g\in G$. In this case, we denote by $G^\Delta$ the permutation group on~$\Delta$, induced by the natural action of~$G$ on~$\Delta$; it consists of the permutations $g^\Delta\in\sym(\Delta)$, $g\in G$, taking~$\delta$ to~$\delta^g$, $\delta\in\Delta$.\medskip

In general, when $G$ acts on a set~$\Delta$ (which is not necessarily a subset of~$\Omega$), we denote by $G^\Delta$ the subgroup of $\sym(\Delta)$, induced by the corresponding action.\medskip

For an arbitrary $\Delta\subseteq\Omega$, we denote by $G_\Delta$ and $G_{\{\Delta\}}$ the pointwise and setwise stabilizer of~$\Delta$ in $G$, respectively. Thus, if $\Delta$ is $G$-invariant, then $G_\Delta$ is the kernel of the restriction homomorphism $G\to G^\Delta$, whereas $G_{\{\Delta\}}=G$. \medskip

The groups $G,H\le\sym(\Delta)$ are said to be $2$-equivalent if $\orb_2(G)=\orb_2(H)$. Formula~\eqref{260920b} implies that $G$ and $\ov G$ are $2$-equivalent. In fact, the $2$-closure $\ov G$ is the largest group in the class of all groups $2$-equivalent to~$G$.\medskip

All other undefined terms and notation are standard and can be found in~\cite{DM}.

\section{Proof of Theorem~\ref{090820a}}\label{sec:Section3}
Let $G\le\sym(\Omega)$ be a nilpotent group of degree $n=|\Omega|$. For a prime~$p\,|\, n$, the largest $p$-power divisor  of~$n$ is denoted by~$n_p$; if $\pi$ is a set of prime divisors of~$n$, then we put $n_\pi:=\prod_{p\in \pi}n_p$. The set of all prime divisors of the order of $G$ is denoted by $\pi(G)$.

\lmml{080820a}
Assume that $G$ is transitive and $H$ a  Hall subgroup of $G$.  Then
\nmrt
\tm{1} the size of every $H$-orbit is equal to $n_\pi$, where $\pi=\pi(H)$,
\tm{2}  $G$ acts on $\orb(H)$; moreover, the kernel of this action is equal to~$H$.
\enmrt
\elmm
\proof The nilpotency of $G$ implies that $H\trianglelefteq G$. Therefore, $G$ permutes the orbits of $H$ and hence acts on~$\Delta=\orb(H)$. Denote by $L$ the kernel of this action.\medskip

The transitivity of $G$ and normality of~$H$ imply that all $H$-orbits have the same size $m$. In particular, $m$ divides~$n$.  Since $n$ divides $|H|$, this implies that  $m$ also  divides~$|H|$. Taking into account that $\GCD(n,|H|)$ divides $n_\pi$,  we conclude that
\qtnl{090820d}
m\ \,\text{divides}\ \, n_\pi.
\eqtn
Next, the group $G^\Delta\le\sym(\Delta)$ is transitive. Consequently, $|\Delta|$ divides~$|G^\Delta|$.  However the numbers $|G^\Delta|$ and $|H|$ are coprime, because $H\le L$ and $H$ is a Hall subgroup of~$G$. Since also~$|\Delta|$ divides~$n$, this implies that $|\Delta|$ divides $n_{\pi'}$, where $\pi'$ is the set of prime divisors of $n$, not belonging to~$\pi$. Together with~\eqref{090820d} this shows that
$$
n=m\cdot|\Delta|\le n_\pi \cdot n_{\pi'}=n,
$$
whence $m=n_\pi$. This proves  statement~(1). \medskip

To prove statement~(2), assume on the contrary that $H<L$. Then there exists an element $g\in L$ of prime order $q\not\in \pi$. Denote by $\Gamma$ the orbit of~$L$ such that $g^\Gamma\ne\id_\Gamma$. Then the Sylow $q$-subgroup~$Q$ of~$L^\Gamma$ is nontrivial. Since $L^\Gamma$ is transitive and nilpotent, we can apply statement~(1) to $G=L^\Gamma$ and $H=Q$. Then
$$
|\Gamma|=|\Gamma|_q\cdot|\orb(Q)|.
$$
On the other hand, we can apply statement~(1) to $G=L^\Gamma$ and $H=H^\Gamma$.  Since $\Gamma$ is also an orbit of~$H$,  this statement implies that
$$
|\Gamma|=|\Gamma|_\pi\cdot|\orb(H^\Gamma)|=|\Gamma|_\pi.
$$
Thus, $|\Gamma|_q\cdot|\orb(Q)|=|\Gamma|_\pi$. However, the numbers $|\Gamma|_q$ and $|\Gamma|_\pi$ are relatively prime. Therefore $|\Gamma|_q=1$ and hence $|\Gamma|=|\orb(Q)|$. It follows that $Q=1$, a contradiction.\eprf

\lmml{100820a}
Theorem~\ref{090820a} holds if $G$ is transitive.
\elmm
\proof If $G$ is a $p$-group, then $\ov G$ is a $p$-group \cite[Exercise~5.28]{WieInvRel} and the required statement is true. Thus we may assume that $G=P\times H$, where $P$ and $H$ are the Sylow $p$-subgroup and Hall subgroup of~$G$, respectively. Let $\Delta\in\orb(P)$ and $\Gamma\in\orb(H)$. By Lemma~\ref{080820a}(1), we  have
$$
|\Delta|=n_{p^{}} \qaq |\Gamma|=n_{p'},
$$
where $p'=\pi(H)$.  Further, $\Delta$ and $\Gamma$ are blocks of the transitive group~$G$. Therefore, the intersection $\Delta\cap\Gamma$ is either empty or is a block the size of which divides both~$|\Delta|$ and $|\Gamma|$. Thus, $|\Delta\cap\Gamma|\le 1$.  \medskip

Each point $\alpha\in\Omega$  lies in exactly one $P$-orbit, say $\Delta_\alpha$,  and in exactly one $H$-orbit, say $\Gamma_\alpha$. By the above, $|\Delta_\alpha\cap\Gamma_\alpha|=1$. Consequently, the mapping
$$
\rho:\Omega\to \orb(P)\times\orb(H),\ \alpha\mapsto(\Delta_\alpha,\Gamma_\alpha)
$$
is a bijection. Denote by $P'$ and $H'$ the permutation groups induced by the actions of $G$ on $\orb(H)$ and on $\orb(P)$, respectively. By  Lemma~\ref{080820a}(2), we have
$$
P^\rho=P'\times 1\qaq H^\rho=1\times H'.
$$
Thus the group $G$ can be identified  with the direct product $P'\times H'$ acting on $\orb(H)\times\orb(P)$. From \cite[Proposition 3.1(2)]{EvdokimovP2001}, it follows that
$$
(P'\times H')^{(2)}=(P')^{(2)}\times (H')^{(2)},
$$
which completes the proof by induction on the number $|\pi(G)|$ with taking into account that $\pi(G)=\pi(P)\cup\pi(H)$.\eprf\medskip

Let us return to the general case. Any transitive constituent $H$ of~$G$ is a homomorphic image of~$G$ which is nilpotent.  Therefore, $H$ is nilpotent and $\pi(H)\subseteq \pi(G)$. By Lemma~\ref{100820a}, the group~$\ov H$ is nilpotent and also $\pi(H)=\pi(\ov H)$. In view of~\cite[Proposition 3.1(1)]{EvdokimovP2001}, this implies that the group
$$
\ov G \le \big( \prod_H H\big)^{(2)}=
\prod_H \ov H
$$
is a subgroup of the direct product of nilpotent groups. Consequently, $\ov G$ is nilpotent and $\pi(G)=\pi(\ov G)$. Thus
it remains to verify that if $P$ and $Q$ are Sylow $p$-subgroups of $G$ and $\ov G$, respectively,  then $\ov P=Q$. We need two auxiliary lemmas.

\lmml{110820l}
Let $P$ and $Q$ as above. Then $\ov P\le Q$.  Moreover,
\qtnl{110820x}
\orb(P)=\orb(Q).
\eqtn
\elmm
\proof The first part of the statement follows from the monotonicity of the $2$-closure operator; in particular, $P\le \ov P\le Q$.  Thus, to prove the second part, it suffices to verify that each $P$-orbit $\Delta$ is a $Q$-orbit. Denote by $\Gamma$ the $G$-orbit containing~$\Delta$. Then $\Gamma$ is also  a $\ov G$-orbit. It follows that there exists a $Q$-orbit $\Delta'$ such that
$$
\Delta\subseteq \Delta'\subseteq\Gamma.
$$
The group $G^\Gamma$ and ${\ov G}\phmb{\Gamma}$ are transitive and nilpotent. By Lemma~\ref{080820a}(1), this implies that $|\Delta|=|\Gamma|_p=|\Delta'|$.  This shows that $\Delta=\Delta'$, as required.\eprf\medskip

\lmml{030619d}
Let $\Delta,\Gamma\in \orb(P)$. Then $(G_{\{\Delta\}}\cap G_{\{\Gamma\}})^{\Delta\cup\Gamma}\le P^{\Delta\cup\Gamma}$.
\elmm
\proof The statement is trivial if $G=P$. Now let $G=P\times H$, where $H$ is a Hall subgroup of~$G$. It follows that  each $g\in G$ can be written as $g=xy$ with $x\in P$ and $y\in H$.  Assume that
$$
g\in G_{\{\Delta\}}\cap G_{\{\Gamma\}},
$$
 i.e.,  $g$ leaves $\Delta$ and $\Gamma$ fixed (as sets). Then the permutation~$x^{-1}\in P$ leaves the sets $\Delta$ and $\Gamma$ fixed, because they are $P$-orbits. Thus, $\Delta^y=\Delta^{x^{-1}g}=\Delta^g=\Delta$ and similarly, $\Gamma^y=\Gamma$.\medskip

We claim that
\qtnl{070820o}
y^\Delta=\id_\Delta\qaq y^\Gamma=\id_\Gamma.
\eqtn
Let us prove the first equality; the second one is  proved analogously. The permutation $y^\Delta$ belongs to the centralizer $Z$ of transitive group $P^\Delta$ in $\sym(\Delta)$. According to~\cite[Exercise~4.5']{Wielandt1964}, the group $Z$ is semiregular. In particular, $|Z|$ divides~$|\Delta|$ which is  a $p$-power. Therefore, $Z$ is a~$p$-group. Consequently, the order of $y^\Delta$  is a $p$-power and hence $y^\Delta\in P^\Delta$. Taking into account that $P^\Delta\cap H^\Delta=1$, we conclude that the first equality in~\eqref{070820o} holds.\medskip

Using equalities \eqref{070820o}, we have
$$
g^{\Delta\cup\Gamma}=(xy)^{\Delta\cup\Gamma}=x^{\Delta\cup\Gamma}y^{\Delta\cup\Gamma}=x^{\Delta\cup\Gamma}\in P^{\Delta\cup\Gamma},
$$
as required.\eprf\medskip

To complete the proof of the theorem,  we note that  $Q$ is $2$-closed, for otherwise~$\ov Q$ is a $p$-subgroup of $\ov G$, containing the Sylow $p$-subgroup~$Q$. Thus it suffices to verify that given $\alpha,\beta\in\Omega$ and $\ov g\in Q$,  there exists $h\in P$ such that
$$
(\alpha,\beta)^{\ov g}=(\alpha,\beta)^h.
$$
The $2$-equivalence of $G$ and $\ov G$ implies that $(\alpha,\beta)^{\ov g}=(\alpha,\beta)^g$ for some $g\in G$.  Denote by $\Delta$ and $\Gamma$ the $Q$-orbits containing the points $\alpha$ and $\beta$, respectively. Then obviously,
$\alpha^g=\alpha^{\ov g}$ belongs to $\Delta$ and  $\beta^g=\beta^{\ov g}$ belongs to~$\Gamma$. Therefore,
$$
\alpha,\alpha^g\in\Delta\qaq \beta,\beta^g\in\Gamma.
$$
In view of equality~\eqref{110820x},  $\Delta$ and $\Gamma$ are also $P$-orbits. Since the group $G$ permute the $P$-orbits, it follows that  $g\in  G_{\{\Delta\}}\cap G_{\{\Gamma\}}$. By Lemma~\ref{030619d}(1), there exists $h\in P$ such that $h^{\Delta\cup\Gamma}=g^{\Delta\cup\Gamma}$. Thus,
$$
(\alpha,\beta)^{\ov g}=(\alpha,\beta)^g=(\alpha,\beta)^h,
$$
as required. Theorem~\ref{090820a} is completely proved.\medskip

{\bf Proof of Corollary~\ref{020619a}.} Let $G$ be a nilpotent permutation group. Assume that~$G$ is $2$-closed. Then $G=\ov G$ and hence $\syl(G)=\syl(\ov G)$. It remains to note that by Theorem~\ref{090820a}, any Sylow subgroup of~$\ov G$ is $2$-closed. Conversely, assume that $P=\ov P$ for each $P\in\syl(G)$. Then  again by Theorem~\ref{090820a}, we have
$$
\ov G=\prod_{P\in\syl(G)}\ov P=\prod_{P\in\syl(G)}P=G,
$$
i.e., $G$ is $2$-closed.\eprf

\section{Proof of Theorem~\ref{200420b}}\label{sec:Section4}

It is well known that a transitive abelian group is regular \cite[Proposition~4.4]{Wielandt1964}. Therefore every abelian permutation group $G$ is {\it quasiregular}, i.e., every transitive constituent of~$G$ is regular. Thus Theorem~\ref{200420b} is an immediate consequence of the following theorem.

\thrml{280820f}
Let $G$ be a quasiregular permutation group and $Z=\zel(G)$. Then~$G$ is $2$-closed if and only if $Z\le G$ and $G^{\orb(Z)}$ is $2$-closed.
\ethrm
\proof Let $\Delta$ be a $G$-orbit. Then the group $G^\Delta$ is regular and hence $2$-closed. Since~$\ov G$ is contained in the direct product of the $2$-closures of the groups  $G^\Delta$, this shows that~$\ov G$ is quasiregular. In what follows, we assume that $G\le\sym(\Omega)$.

\lmml{280820g}
$Z\le \ov G$, and also $Z^\Delta\trianglelefteq \ov G\phmb{\Delta}$ for all $\Delta\in\orb(\ov G)$.
\elmm
\proof The group $Z$ is generated by the permutations~$z$ satisfying the following condition: there exists $\Delta\in\orb(Z)$ such that
\qtnl{290820a}
z^\Delta\in Z^\Delta\qaq z^{\Omega\setminus \Delta}=\id_{\Omega\setminus\Delta}.
\eqtn
Thus to prove the inclusion $Z\le\ov G$, it suffices to verify that each such $z$ belongs to~$\ov G$, or equivalently that $s^z=s$ for every $s\in\orb_2(G)$. \medskip

Denote by $\Delta'$ and $\Delta''$ the $G$-orbits such that $s\subseteq\Delta'\times\Delta''$. If $\Delta'\ne \Delta\ne\Delta''$, then $s^z=s$ by the second equality in~\eqref{290820a}, whereas if $\Delta=\Delta'=\Delta''$, then $s\in\orb_2(G^\Delta)$ and again $s^z=s$, because $Z^\Delta\le G^\Delta$. Thus without loss of generality, we may assume that $\Delta=\Delta''\ne\Delta'$. Then by the definition of $Z$ (see formula~\eqref{270820a}), there exists $g\in G_{\Delta'}$ such that
$$
z^{\Delta\cup\Delta'}=g^{\Delta\cup\Delta'}.
$$
Thus, $s^z=s^g=s$, as required. \medskip

Let us prove that $Z^\Delta\trianglelefteq \ov G\phmb{\Delta}$. We have $G_{\Delta'}\trianglelefteq G$ for every  $\Delta'\in\orb(G)$. Therefore, $G_{\Delta'}^\Delta\trianglelefteq G^\Delta$. Consequently, the intersection of all $G_{\Delta'}^\Delta$ taken over all $\Delta'\ne\Delta$ is also normal in~$G^\Delta$. Since this intersection coincides with $Z^\Delta$ and $G^\Delta=\ov G\phmb{\Delta}$, we are done.\eprf\medskip

From Lemma~\ref{280820g}, it follows that for every $\Delta\in\orb(\ov G)$, the group $\ov G\phmb{\Delta}$  acts on the set $\orb(Z^\Delta)$. Therefore, $\ov G$ acts on the union $\orb(Z)$ of all $\orb(Z^\Delta)$. Denote by~$\ov\rho$ the corresponding epimorphism from~$\ov G$ to~$\ov G\phmb{\orb(Z)}$. Then  obviously~$Z$ is a subgroup of $L:=\ker(\ov\rho)$.  Moreover, if $\Delta$ is a $\ov G$-orbit, then
$$
L^\Delta=Z^\Delta,
$$
because $L^\Delta$ and $Z^\Delta$ are $1$-equivalent subgroups of the regular group $\ov G\phmb{\Delta}$ (recall that the group~$\ov G$ is quasiregular). Since $L$ is contained in the direct product of the~$Z^\Delta$, the definition of $Z$ implies that $L=Z$. This proves the following statement.

\lmml{280820a}
$\ker(\ov\rho)=Z$.
\elmm

The epimorphism $\ov\rho$ induces the action of $G\le \ov G$ on the set $\orb(Z)$; denote by~$\rho$ the corresponding epimorphism from $G$ to~$G\phmb{\orb(Z)}$. It should be noted that while the group $Z$ is not, in general, a subgroup of~$G$, the permutation group $G\phmb{\orb(Z)}$ is well defined.

\lmml{280820b}
$\im(\ov\rho)=\ov{G\phmb{\orb(Z)}}$.
\elmm
\proof The groups $G\phmb{\orb(Z)}$ and $\ov G\phmb{\orb(Z)}$ are $2$-equivalent \cite[Lemma~2.1(1)]{VasP2020}. Therefore,
$$
\im(\ov\rho)=\ov G\phmb{\orb(Z)}\le \ov{G\phmb{\orb(Z)}}.
$$
Conversely, we need to verify that for every $\ov g\in \ov{G\phmb{\orb(Z)}}$, there exists $g\in\ov G$ such that
\qtnl{300620a}
\ov\rho(g)=\ov g.
\eqtn

Let $\Delta\in\orb(\ov G)$. The quasiregularity of~$\ov G$ implies that the group $\ov G\phmb{\Delta}=G^\Delta$ is regular. Since also  $Z^\Delta\trianglelefteq \ov G\phmb{\Delta}$ (Lemma~\ref{280820g}), the group $(\ov G\phmb{\Delta})\phmb{\orb(Z^\Delta)}$ is also regular and hence $2$-closed. It follows that
$$
\bigl(\ov{G\phmb{\orb(Z)}}\bigr)\phmb{\ov\Delta}\le
\ov{(G\phmb{\Delta})\phmb{\orb(Z^\Delta)}}=
\ov{(\ov{G}\phmb{\Delta})\phmb{\orb(Z^\Delta)}}=
(\ov{G}\phmb{\Delta})\phmb{\orb(Z^\Delta)},
$$
where $\ov\Delta$ is the $\ov{G\phmb{\orb(Z)}}$-orbit the points of which are the $Z$-orbits contained in~$\Delta$. Thus for every $\Delta\in\orb(\ov G)$, there exists a permutation $g_\Delta\in \ov G\phmb{\Delta}$ such that
$$
\ov\rho_\Delta(g_\Delta)={\ov g}\phmb{\,\ov\Delta},
$$
where $\ov\rho_\Delta$ is the epimorphism from $\ov G\phmb{\Delta}$ to $(\ov G\phmb{\Delta})\phmb{\orb(Z^\Delta)}$,  induced by~$\ov\rho$.  Now if the  product~$g$ of all  the $g_\Delta$ lies in~$G$, then
$$
\ov\rho(g)=\ov\rho\bigl(\prod_\Delta g_\Delta\bigr)=\prod_{\Delta} \ov\rho_\Delta(g_\Delta)=
\prod_{\ov\Delta} \ov g\phmb{\ov\Delta}=\ov g,
$$
which proves equality~\eqref{300620a}.\medskip

It remains to verify that $g\in\ov G$. To this end, let $s\in\orb_2(G)$. Then $s\in\orb_2(\ov G)$ and hence $s\subseteq \Delta\times\Gamma$ for some $\Delta,\Gamma\in\orb(\ov G)$. Now if $\Delta=\Gamma$, then $s^g=s^{g_\Delta}=s$, because $g_\Delta\in\ov G\phmb{\Delta}$. Assume that $\Delta\ne\Gamma$. Then by Lemma~\ref{280820a} and the definition of~$Z$, we have $Z^\Delta\times Z^\Gamma\subseteq \ov G\phmb{\Delta\cup\Gamma}$. It follows that
$$
(\alpha,\beta)\in s\quad\Leftrightarrow\quad \ov\alpha\times\ov\beta\subseteq s,
$$
where $\ov\alpha=\alpha^Z$ and $\ov\beta=\beta^Z$. Furthermore, the set $\ov s=\{(\ov\alpha,\ov\beta):\ (\alpha,\beta)\in s\}$ is a $2$-orbit of the group $\ov G\phmb{\orb(Z)}$ and hence $\ov s\phmb{\ov g}=\ov s$. Thus,
$$
s^g=\bigl(\bigcup_{(\ov\alpha,\ov\beta)\in\ov s}\ov\alpha\times\ov\beta\,\bigr)^g=
\bigcup_{(\ov\alpha,\ov\beta)\in\ov s}\ov\alpha\phmb{\ov g}\times\ov\beta\phmb{\ov g}=
\bigcup_{(\ov\alpha,\ov\beta)\in\ov s\phmb{\,\ov g}}\ov\alpha\times\ov\beta=
\bigcup_{(\ov\alpha,\ov\beta)\in\ov s}\ov\alpha\times\ov\beta=s
$$
as required. \eprf\medskip

To prove the ``only if'' part, assume that the group $G$ is $2$-closed. Then by Lemma~\ref{280820g}, we have $Z\le \ov G=G$, whereas by Lemma~\ref{280820b}, we have
$$
G\phmb{\orb(Z)}=\ov G\phmb{\orb(Z)}=\im(\ov\rho)=\ov{G\phmb{\orb(Z)}},
$$
i.e., the group $G\phmb{\orb(Z)}$ is $2$-closed, as required.\medskip

To prove the ``if'' part, assume that $Z\le G$ and the group $G\phmb{\orb(Z)}$ is $2$-closed.  By Lemma~\ref{280820a}, the first condition implies that
$$
Z\le\ker(\rho)\le\ker(\ov\rho)=Z,
$$
in particular, $\ker(\rho)=\ker(\ov\rho)$. Furthermore, by Lemma~\ref{280820b} the second condition implies that $\ov{G\phmb{\orb(Z)}}=G\phmb{\orb(Z)}$ and hence
$$
\im(\rho)=G\phmb{\orb(Z)}=\ov{G\phmb{\orb(Z)}}=\im(\ov\rho).
$$
Thus,
$$
|G|=|\ker(\rho)|\cdot|\im(\rho)|=|\ker(\ov\rho)|\cdot|\im(\ov\rho)|=|\ov G|.
$$
Since $G\le\ov G$, this means that $G=\ov G$, i.e., $G$ is $2$-closed.\eprf

\section{Proof of Theorem~\ref{310820a}}\label{260920a}

We begin with a sufficient condition for an orbit of quasiregular permutation group to be unessential. The proof is based on a special result from theory of coherent configurations~\cite{CCBook}.

\lmml{220920a}
Let $G$ be a  quasiregular permutation group and $\Delta\in \orb(G)$. Assume that $G^\Delta_{\Delta'}=1$ for some $G$-orbit $\Delta'\ne\Delta$. Then the orbit~$\Delta$ is unessential.
\elmm
\proof The quasiregularity of~$G$ implies that  $G_{\delta'}=G_{\Delta'}$ for each $\delta'\in\Delta'$. It follows that if $\delta,\lambda\in\Delta$, then
\qtnl{210920a}
(\delta',\delta)\in (\delta',\lambda)^G\quad\Rightarrow\quad \delta=\lambda.
\eqtn
Indeed, if $(\delta',\delta)=(\delta',\lambda)^g$ for some $g\in G$, then $g\in G_{\delta'}=G_{\Delta'}$. Since $G^\Delta_{\Delta'}=1$, this implies that $\delta=\lambda^g=\lambda$.\medskip

Denote by $\Omega$ the point set of~$G$  and put $S=\orb_2(G)$. Then the pair $\cX=(\Omega,S)$ is a coherent configuration and $\ov G=\aut(\cX)$ is the automorphism group of~$\cX$, see~\cite{CCBook}. Formula~\eqref{210920a} implies that the condition~(3.3.14) from~\cite{CCBook} is satisfied for the coherent configuration~$\cX$ and the set $\Delta$ equal to
$$
\Omega':=\Omega\setminus\Delta.
$$
By~\cite[Lemma~3.3.20(1)]{CCBook}, the restriction homomorphism $\aut(\cX)\to\aut(\cX_{\Omega'})$ is an isomorphism; in particular,
\qtnl{210920x}
\ov G\phmb{\Omega'}=\aut(\cX)^{\Omega'}=\aut(\cX_{\Omega'})=\ov{G^{\Omega'}}.
\eqtn

To prove that the orbit $\Delta$ is unessential, first assume that the group $G$ is $2$-closed. Then formula~\eqref{210920x} shows that the group $G^{\Omega'}=\ov G\phmb{\Omega'}$ is also $2$-closed. Conversely, assume that $G^{\Omega'}$ is $2$-closed. Then $G^{\Omega'}=\ov{G^{\Omega'}}$. Consequently,  the restriction homomorphism $G\to G^{\Omega'}$ is surjective. Since $G\le \ov G=\aut(\cX)$, the homomorphism is also injective. Thus,  $G=\ov G$, i.e., $G$ is $2$-closed. \eprf\medskip

Turn to the proof of Theorem~\ref{310820a}. Assume that the group $\zel(G)$ is trivial. Then for each $\Delta\in\orb(G)$, we have
$$
\bigcap\limits_{\Delta'\ne\Delta}G^\Delta_{\Delta'}=1.
$$
On the other hand, by the theorem hypothesis, $G^\Delta_{\Delta'}\le G^\Delta_{}$ is a cyclic $p$-group for all $G$-orbits~$\Delta'$. Thus, $G^\Delta_{\Delta'}=1$ for at least one $\Delta'$. By Lemma~\ref{220920a}, this implies that  the orbit~$\Delta$ is unessential.

\bigskip
\textsc{Dmitry Churikov\\
Sobolev Institute of Mathematics,\\
pr. Koptyuga, 4,\\
630090, Novosibirsk, Russia\\
Novosibirsk State University,\\
ul. Pirogova, 2,\\
630090, Novosibirsk, Russia\\}
\verb"churikovdv@gmail.com"\\

\noindent\textsc{Ilia Ponomarenko\\
Steklov Institute of Mathematics at St. Petersburg,\\
Fontanka, 27,\\
191023, St. Petersburg, Russia\\
Sobolev Institute of Mathematics,\\
pr. Koptyuga, 4,\\
630090, Novosibirsk, Russia\\}
\verb"inp@pdmi.ras.ru"\\
\end{document}